\documentclass[12pt,a4paper]{article}
\usepackage{amsmath}
\usepackage{amssymb}
\usepackage{t1enc}
\usepackage[latin1]{inputenc}
\usepackage[english]{babel}
\usepackage{amsfonts}
\usepackage{latexsym}
\usepackage{bm}
\newtheorem{theorem}{Theorem}[section]
\newtheorem{prop}[theorem]{Proposition}
\newtheorem{lemma}[theorem]{Lemma}
\newtheorem{cor}[theorem]{Corollary}
\newtheorem{conj}[theorem]{Conjecture}

\begin{document}
\title{Cross-intersecting sub-families of hereditary families}

\author{Peter Borg\\[5mm]
Department of Mathematics, University of Malta, Msida MSD 2080,
Malta\\
\texttt{p.borg.02@cantab.net}}
\date{\today} \maketitle

\begin{abstract}
Families $\mathcal{A}_1, \mathcal{A}_2, ..., \mathcal{A}_k$ of
sets are said to be \emph{cross-intersecting} if for any $i$ and
$j$ in $\{1, 2, ..., k\}$ with $i \neq j$, any set in
$\mathcal{A}_i$ intersects any set in $\mathcal{A}_j$. For a
finite set $X$, let $2^X$ denote the \emph{power set of $X$} (the
family of all subsets of $X$). A family $\mathcal{H}$ is said to
be \emph{hereditary} if all subsets of any set in $\mathcal{H}$
are in $\mathcal{H}$; so $\mathcal{H}$ is hereditary if and only
if it is a union of power sets. We conjecture that for any
non-empty hereditary sub-family $\mathcal{H} \neq \{\emptyset\}$
of $2^X$ and any $k \geq |X|+1$, both the sum and product of
sizes of $k$ cross-intersecting sub-families $\mathcal{A}_1,
\mathcal{A}_2, ..., \mathcal{A}_k$ (not necessarily distinct or
non-empty) of $\mathcal{H}$ are maxima if $\mathcal{A}_1 =
\mathcal{A}_2 = ... = \mathcal{A}_k = \mathcal{S}$ for some
largest \emph{star $\mathcal{S}$ of $\mathcal{H}$} (a sub-family
of $\mathcal{H}$ whose sets have a common element).
We prove this for the case when $\mathcal{H}$ is \emph{compressed
with respect to an element $x$ of $X$}, 
and for this purpose we establish new properties of the usual
\emph{compression operation}. For the product, we actually
conjecture that the configuration $\mathcal{A}_1 = \mathcal{A}_2
= ... = \mathcal{A}_k = \mathcal{S}$ is optimal for any
hereditary $\mathcal{H}$ and any $k \geq 2$, and we prove this
for a special case too.
\end{abstract}

\section{Basic definitions and notation} \label{Intro}
Unless otherwise stated, we shall use small letters such as $x$ to
denote elements of a set or non-negative integers or functions,
capital letters such as $X$ to denote sets, and calligraphic
letters such as $\mathcal{F}$ to denote \emph{families}
(i.e.~sets whose elements are sets themselves). It is to be
assumed that sets and families are \emph{finite}. We call a set
$A$ an \emph{$r$-element set}, or simply an \emph{$r$-set}, if
its size $|A|$ is $r$ (i.e.~if it contains exactly $r$ elements).

For any integer $n \geq 1$, the set $\{1, ..., n\}$ of the first
$n$ positive integers is denoted by $[n]$. For a set $X$, the
\emph{power set of $X$} (i.e.~$\{A \colon A \subseteq X\}$) is
denoted by $2^X$, and the family of all $r$-element subsets of
$X$ is denoted by $X \choose r$.

A family $\mathcal{H}$ is said to be a \emph{hereditary family}
(also called an \emph{ideal} or a \emph{downset}) if all the
subsets of any set in $\mathcal{H}$ are in $\mathcal{H}$. Clearly
a family is hereditary if and only if it is a union of power sets.
A \emph{base of $\mathcal{H}$} is a set in $\mathcal{H}$ that is
not a subset of any other set in $\mathcal{H}$. So a hereditary
family is the union of power sets of its bases. An interesting
example of a hereditary family is the family of all independent
sets of a graph or matroid.

We will denote the union of all sets in a family $\mathcal{F}$ by
$U(\mathcal{F})$. If $x$ is an element of a set $X$, then we
denote the family of those sets in $\mathcal{F}$ which contain
$x$ by $\mathcal{F} \langle x \rangle$, and we call $\mathcal{F}
\langle x \rangle$ a \emph{star of $\mathcal{F}$}. So $\mathcal{F}
\langle x \rangle$ is the empty set $\emptyset$ if and only if
$x$ is not in $U(\mathcal{F})$.


A family $\mathcal{A}$ is said to be \emph{intersecting} if any
two sets 
in $\mathcal{A}$ intersect (i.e.~contain at least one common
element). We call a family $\mathcal{A}$ \emph{centred} if the
sets in $\mathcal{A}$ have a common element $x$
(i.e.~$\mathcal{A} = \mathcal{A}\langle x \rangle$). So a centred
family is intersecting, and a non-empty star of a family
$\mathcal{F}$ is centred. The simplest example of a non-centred
intersecting family is $\{\{1,2\}, \{1,3\}, \{2,3\}\}$
(i.e.~${[3] \choose 2}$).

Families $\mathcal{A}_1, ..., \mathcal{A}_k$ are said to be
\emph{cross-intersecting} if for any $i$ and $j$ in $[k]$ with $i
\neq j$, any set in $\mathcal{A}_i$ intersects any set in
$\mathcal{A}_j$.

If $U(\mathcal{F})$ has an element $x$ such that $\mathcal{F}
\langle x \rangle$ is a largest intersecting sub-family of
$\mathcal{F}$ (i.e.~no intersecting sub-family of $\mathcal{F}$
has more sets than $\mathcal{F}\langle x \rangle$), then we say
that $\mathcal{F}$ has the \emph{star property at $x$}. We simply
say that $\mathcal{F}$ has the \emph{star property} if either
$U(\mathcal{F}) = \emptyset$ or $\mathcal{F}$ has the star
property at some element of $U(\mathcal{F})$.

If $U(\mathcal{F})$ has an element $x$ such that $(F \backslash
\{y\}) \cup \{x\} \in \mathcal{F}$ whenever $y \in F \in
\mathcal{F}$ and $x \notin F$, then $\mathcal{F}$ is said to be
\emph{compressed with respect to $x$}. For example, this is the
case when $\mathcal{F}$ is the family of all independent sets of a
graph that has an isolated vertex $x$.

A family $\mathcal{F} \subseteq 2^{[n]}$ is said to be
\emph{left-compressed} if $(F \backslash \{j\}) \cup \{i\} \in
\mathcal{F}$ whenever $1 \leq i < j \in F \in \mathcal{F}$ and $i
\notin F$.

\section{Intersecting sub-families of hereditary families}
\label{Her}

The following is a famous longstanding open conjecture in
extremal set theory due to Chv\'{a}tal (see \cite{Borg3} for a
more general conjecture).

\begin{conj}[\cite{Chv}] \label{Chvatal} If $\mathcal{H}$
is a hereditary family, then $\mathcal{H}$ has the star property.
\end{conj}
This conjecture was verified for the case when $\mathcal{H}$ is
left-compressed by Chv\'{a}tal \cite{Chva} himself. Snevily
\cite{Sn} took this result (together with results in
\cite{Schonheim, Wang}) a significant step forward by verifying
Conjecture~\ref{Chvatal} for the case when $\mathcal{H}$ is
compressed with respect to an element $x$ of $U(\mathcal{H})$.

\begin{theorem} [\cite{Sn}] \label{Snevily} If a hereditary family
$\mathcal{H}$ is compressed with respect to an element $x$ of
$U(\mathcal{H})$, then $\mathcal{H}$ has the star property at $x$.
\end{theorem}
%
A generalisation is proved in \cite{Borg3} by means of an
alternative self-contained argument.

Snevily's proof of Theorem~\ref{Snevily} makes use of the
following interesting result of Berge \cite{Berge} (a proof of
which is also provided in \cite[Chapter 6]{Anderson}).
\begin{theorem} [\cite{Berge}] \label{Berge} If
$\mathcal{H}$ is a hereditary family, then $\mathcal{H}$ is a
disjoint union of pairs of disjoint sets, together with
$\emptyset$ if $|\mathcal{H}|$ is odd.
\end{theorem}
This result was also motivated by Conjecture~\ref{Chvatal} as it
has the following immediate consequence.

\begin{cor}\label{Bergecor} If $\mathcal{A}$ is an intersecting
sub-family of a hereditary family $\mathcal{H}$, then
$$|\mathcal{A}| \leq \frac{1}{2}|\mathcal{H}|.$$
\end{cor}
\textbf{Proof.} For any pair of disjoint sets, at most only one
set can be in an intersecting family $\mathcal{A}$. By
Theorem~\ref{Berge}, the result follows.~\hfill{$\Box$} \\

A special case of Theorem~\ref{Snevily} is a result of Schönheim
\cite{Schonheim} which says that Conjecture~\ref{Chvatal} is true
if the bases of $\mathcal{H}$ have a common element, and this
follows immediately from Corollary~\ref{Bergecor} and the
following fact.

\begin{prop} \label{Bergeprop} If the bases of a hereditary family
$\mathcal{H}$ have a common element $x$, then $$|\mathcal{H}
\langle x \rangle| = \frac{1}{2}|\mathcal{H}|.$$
\end{prop}
\textbf{Proof.} By induction on $|U(\mathcal{H})|$. We have $x
\in U(\mathcal{H})$. Let $\mathcal{A} = \mathcal{H}\langle x
\rangle$. If $|U(\mathcal{H})| = 1$, then $\mathcal{H} =
\{\emptyset, \{x\}\}$, $\mathcal{A} = \{\{x\}\}$ and hence
$|\mathcal{A}| = 1 = \frac{1}{2}|\mathcal{H}|$. Now suppose
$|U(\mathcal{H})| \geq 2$. Then $U(\mathcal{H})$ has an element $y
\neq x$. Let $\mathcal{I} = \{H \backslash \{y\} \colon H \in
\mathcal{H}\langle y \rangle\}$ and $\mathcal{J} = \{H \in
\mathcal{H} \colon y \notin H\}$. Similarly, let $\mathcal{B} =
\{A \backslash \{y\} \colon A \in \mathcal{A}\langle y \rangle\}$
and $\mathcal{C} = \{A \in \mathcal{A} \colon y \notin A\}$. Note
that $\mathcal{I}$ and $\mathcal{J}$ are hereditary, the bases of
$\mathcal{I}$ and $\mathcal{J}$ contain $x$, $|U(\mathcal{I})|
\leq |U(\mathcal{H})| - 1$, $|U(\mathcal{J})| \leq
|U(\mathcal{H})| - 1$, $\mathcal{B} = \mathcal{I} \langle x
\rangle$ and $\mathcal{C} = \mathcal{J} \langle x \rangle$. By
the inductive hypothesis, $|\mathcal{B}| =
\frac{1}{2}|\mathcal{I}|$ and $|\mathcal{C}| =
\frac{1}{2}|\mathcal{J}|$. Finally, $|\mathcal{A}| = |\mathcal{A}
\langle y \rangle| + |\mathcal{C}| = |\mathcal{B}| +
|\mathcal{C}| = \frac{1}{2}(|\mathcal{I}| + |\mathcal{J}|) =
\frac{1}{2}(|\mathcal{H} \langle y \rangle| + |\mathcal{J}|) =
\frac{1}{2}|\mathcal{H}|$.~\hfill{$\Box$}\\

Many other results and problems have been inspired by
Conjecture~\ref{Chvatal} or are related to it; see
\cite{Chvatalsite, Miklos, West}.

\section{Cross-intersecting sub-families of hereditary families}

For intersecting sub-families of a given family $\mathcal{F}$, the
natural question to ask is how large they can be.
Conjecture~\ref{Chvatal} claims that when $\mathcal{F}$ is
hereditary, we need only check the stars of $\mathcal{F}$ (of
which there are $|U(\mathcal{F})|$). For cross-intersecting
families, two natural parameters arise: the sum and the product
of sizes of the cross-intersecting families (note that the product
of sizes of $k$ families $\mathcal{A}_1, ..., \mathcal{A}_k$ is
the number of $k$-tuples $(A_1, ..., A_k)$ such that $A_i \in
\mathcal{A}_i$ for each $i \in [k]$). It is therefore natural to
consider the problem of maximising the sum or the product of sizes
of $k$ cross-intersecting sub-families (not necessarily distinct
or non-empty) of a given family $\mathcal{F}$ (see \cite{Borg2}).
We suggest a few conjectures for the case when $\mathcal{F}$ is
hereditary, and we prove that they are true in some important
cases. Obviously, any family $\mathcal{F}$ is a sub-family of
$2^{X}$ with $X = U(\mathcal{F})$, and we may assume that $X =
[n]$.

For the sum of sizes, we suggest the following.

\begin{conj} \label{weaksum} If $k \geq n+1$ and $\mathcal{A}_1, ...,
\mathcal{A}_k$ are cross-intersecting sub-families of a
hereditary sub-family $\mathcal{H} \neq \{\emptyset\}$ of
$2^{[n]}$, then the sum $\sum_{i=1}^k|\mathcal{A}_i|$ is maximum
if $\mathcal{A}_1 = ... = \mathcal{A}_k = \mathcal{S}$ for some
largest star $\mathcal{S}$ of $\mathcal{H}$.
\end{conj}
We cannot remove the condition that $k \geq n+1$. Indeed,
consider $\mathcal{H} = \{\emptyset\} \cup {[n] \choose 1}$ and
$2 \leq k < n+1$. Let $\mathcal{S} = \{\{1\}\}$; so $\mathcal{S}$
is a largest star of $\mathcal{H}$. Let $\mathcal{A}_1 = ... =
\mathcal{A}_k = \mathcal{S}$, and let $\mathcal{B}_1 =
\mathcal{H}$ and $\mathcal{B}_2 = ... = \mathcal{B}_k =
\emptyset$. Then $\mathcal{A}_1, ..., \mathcal{A}_k$ are
cross-intersecting, $\mathcal{B}_1, ..., \mathcal{B}_k$ are
cross-intersecting, and $\sum_{i = 1}^k |\mathcal{A}_i| = k < n+1
= \sum_{i = 1}^k |\mathcal{B}_i|$. Also, we cannot remove the
condition that $\mathcal{H} \neq \{\emptyset\}$. Indeed, suppose
$\mathcal{H} = \{\emptyset\}$; so $\mathcal{S} = \emptyset$.
Thus, if $\mathcal{A}_1 = ... = \mathcal{A}_k = \mathcal{S}$,
$\mathcal{B}_1 = \mathcal{H}$ and $\mathcal{B}_2 = ... =
\mathcal{B}_k = \emptyset$, then $\sum_{i = 1}^k |\mathcal{A}_i|
= 0 < 1 = \sum_{i = 1}^k |\mathcal{B}_i|$.

For the general case when we have any number of
cross-intersecting families, we suggest the following stronger
conjecture.

\begin{conj} \label{strongsum} Let $\mathcal{A}_1, ...,
\mathcal{A}_k$ be cross-intersecting sub-families of a non-empty
hereditary sub-family $\mathcal{H} \neq \{\emptyset\}$ of
$2^{[n]}$, and let $\mathcal{S}$ be a largest star of
$\mathcal{H}$. \\
(i) If $k \leq \frac{|\mathcal{H}|}{|\mathcal{S}|}$, then
$\sum_{i=1}^k|\mathcal{A}_i|$ is maximum if $\mathcal{A}_1 =
\mathcal{H}$ and $\mathcal{A}_2 = ... =
\mathcal{A}_k = \emptyset$. \\
(ii) If $k \geq \frac{|\mathcal{H}|}{|\mathcal{S}|}$, then
$\sum_{i=1}^k|\mathcal{A}_i|$ is maximum if $\mathcal{A}_1 = ...
= \mathcal{A}_k = \mathcal{S}$.
\end{conj}
This conjecture is simply saying that at least one of the two
simple configurations $\mathcal{A}_1 = \mathcal{H}$,
$\mathcal{A}_2 = ... = \mathcal{A}_k = \emptyset$ and
$\mathcal{A}_1 = ... = \mathcal{A}_k = \mathcal{S}$ gives a
maximum sum of sizes. This strengthens Conjecture~\ref{weaksum}
because, since $\mathcal{H}$ has a non-empty set (as $\mathcal{H}
\neq \emptyset$ and $\mathcal{H} \neq \{\emptyset\}$), we have
$\mathcal{S} \neq \emptyset$, $|\mathcal{H}| =
\left|\{\emptyset\} \cup \bigcup_{i=1}^n \mathcal{H}\langle i
\rangle \right| \leq 1+n|\mathcal{S}| \leq (n+1)|\mathcal{S}|$
and hence $\frac{|\mathcal{H}|}{|\mathcal{S}|} \leq n+1$; that
is, if (ii) is true, then Conjecture~\ref{weaksum} follows.

For the product of sizes, we first present the following
consequence of Conjecture~\ref{weaksum}.

\begin{prop} \label{sumprod} Let $A_1, ..., \mathcal{A}_k,
\mathcal{H}$ and $\mathcal{S}$ be as in Conjecture~\ref{weaksum}.
If Conjecture~\ref{weaksum} is true, then the product
$\prod_{i=1}^k|\mathcal{A}_i|$ is maximum if $\mathcal{A}_1 = ...
= \mathcal{A}_k = \mathcal{S}$.
\end{prop}
This follows immediately from the following elementary result,
known as the Arithmetic Mean-Geometric Mean (AM-GM) Inequality.

\begin{lemma}[AM-GM Inequality] \label{AMGM} If $x_1, x_2, ...,
x_k$ are non-negative real numbers, then
\[\left( \prod_{i=1}^k x_i \right)^{1/k} \leq \frac{1}{k}\sum_{i=1}^k
x_i.\]
\end{lemma}
Indeed, suppose Conjecture~\ref{weaksum} is true. Then
$\sum_{i=1}^k|\mathcal{A}_i| \leq k|\mathcal{S}|$. Thus, by
Lemma~\ref{AMGM}, $\left( \prod_{i=1}^k|\mathcal{A}_i|
\right)^{1/k} \leq |\mathcal{S}|$ and hence
Proposition~\ref{sumprod}.


However, we conjecture the following stronger statement about the
maximum product.

\begin{conj} \label{prodconj} If $k \geq 2$ and $\mathcal{A}_1,
..., \mathcal{A}_k$ are cross-intersecting sub-families of a
hereditary family $\mathcal{H}$, then
$\prod_{i=1}^k|\mathcal{A}_i|$ is maximum if $\mathcal{A}_1 = ...
= \mathcal{A}_k = \mathcal{S}$ for some largest star
$\mathcal{S}$ of $\mathcal{H}$.
\end{conj}
If the above conjecture is true for $k = 2$, then it is true for
any $k \geq 2$. Indeed, it is not difficult to show that, in
general, if $p \geq 2$ and $\mathcal{L}$ is an intersecting
sub-family of a family $\mathcal{F}$ such that the product of
sizes of $p$ cross-intersecting sub-families of $\mathcal{F}$ is a
maximum when each of them is $\mathcal{L}$, then for any $k \geq
p$, the product of sizes of $k$ cross-intersecting sub-families
of $\mathcal{F}$ is also maximum when each of them is
$\mathcal{L}$; see \cite{Borg2}.

Each of the above conjectures generalises
Conjecture~\ref{Chvatal}. Indeed, let $\mathcal{A}$ be an
intersecting sub-family of a hereditary family $\mathcal{H}
\subseteq 2^{[n]}$ with $U(\mathcal{H}) \neq \emptyset$, and let
$\mathcal{S}$ be a largest star of $\mathcal{H}$. Let $k \geq
n+1$, and let $\mathcal{A}_1 = ... = \mathcal{A}_k =
\mathcal{A}$. Then $\mathcal{A}_1, ..., \mathcal{A}_k$ are
cross-intersecting. Thus, any of Conjectures~\ref{weaksum},
\ref{strongsum} and \ref{prodconj} claims that $|\mathcal{A}_i|
\leq |\mathcal{S}|$ for each $i \in [k]$ (since $\mathcal{A}_1 =
... = \mathcal{A}_k$), and hence $|\mathcal{A}| \leq
|\mathcal{S}|$ as claimed by Conjecture~\ref{Chvatal}.

All the above conjectures are true for the special case when
$\mathcal{H} = 2^{[n]}$; more precisely, the following holds.

\begin{theorem} [\cite{Borg2}] \label{powersetcase} For any $k
\geq 2$, both the sum and product of sizes of $k$
cross-intersecting sub-families $\mathcal{A}_1, ...,
\mathcal{A}_k$ of $2^{[n]}$ are maxima if $\mathcal{A}_1 = ... =
\mathcal{A}_k = \{A \subseteq [n] \colon 1 \in A\}$.
\end{theorem}
We generalise this result as follows.
\begin{theorem} \label{result2} If $\mathcal{A}_1, ...,
\mathcal{A}_k$ are cross-intersecting sub-families of a hereditary
family $\mathcal{H}$, then
\[\sum_{i = 1}^k |\mathcal{A}_i| \leq k \frac{|\mathcal{H}|}{2}
\quad \mbox{ and } \quad \prod_{i = 1}^k |\mathcal{A}_i| \leq
\left ( \frac{|\mathcal{H}|}{2}\right ) ^k.\]
Moreover, both bounds are attained if the bases of $\mathcal{H}$
have a common element $x$ and $\mathcal{A}_1 = ... =
\mathcal{A}_k = \mathcal{H}\langle x \rangle$.
\end{theorem}
\textbf{Proof.} Theorem~\ref{Berge} tells us that there exists a
partition $\mathcal{H}_1 \cup \mathcal{H}_2 \cup ... \cup
\mathcal{H}_m$ of $\mathcal{H}$ such that $m = \left \lceil
\frac{|\mathcal{H}|}{2} \right \rceil$, $\mathcal{H}_i =
\{H_{i,1}, H_{i,2}\}$ for some $H_{i,1}, H_{i,2} \in \mathcal{H}$
with $H_{i,1} \cap H_{i,2} = \emptyset$, $i = 1, ..., m$, and if
$|\mathcal{H}|$ is odd then $H_{m,1} = H_{m,2} = \emptyset$.

Let $\mathcal{A} = \bigcup_{i = 1}^k \mathcal{A}_i$. By the
cross-intersection condition, we clearly have $\mathcal{A}^* =
\bigcup_{i=1}^k \mathcal{A}_{i}^*$ and $\mathcal{A}' = \bigcup_{i
= 1}^k \mathcal{A}_i'$. Suppose $\mathcal{A}_i' \cap
\mathcal{A}_j' \neq \emptyset$ for some $i \neq j$. Let $A \in
\mathcal{A}_i' \cap \mathcal{A}_j'$.  Then there exists $A_i \in
\mathcal{A}_i'$ such that $A \cap A_i = \emptyset$, but this is a
contradiction because $A \in \mathcal{A}_j$. So $\mathcal{A}_i'
\cap \mathcal{A}_j' = \emptyset$ for any $i \neq j$. Therefore
$|\mathcal{A}'| = \sum_{i=1}^k |\mathcal{A}_i'|$.

%
Let $\mathcal{B} = \{H_{i,j} \colon i \in [m], j \in [2],
H_{i,3-j} \in \mathcal{A}^*\}$. So $|\mathcal{B}| =
|\mathcal{A}^*|$. For any $H_{i,j} \in \mathcal{B}$, $H_{i,j}
\notin \mathcal{A}$ since $H_{i,3-j} \in \mathcal{A}^*$ and
$H_{i,j} \cap H_{i,3-j} = \emptyset$. So $\mathcal{A}$ and
$\mathcal{B}$ are disjoint sub-families of $\mathcal{H}$.
Therefore,
\begin{equation} 2|\mathcal{A}^*| + |\mathcal{A}'| =
|\mathcal{A}^*| + |\mathcal{B}| + |\mathcal{A}'| = |\mathcal{A}|
+ |\mathcal{B}| = |\mathcal{A} \cup \mathcal{B}| \leq
|\mathcal{H}| \nonumber
\end{equation}
and hence, dividing throughout by $2$, we get $|\mathcal{A}^*| +
\frac{1}{2} |\mathcal{A}'| \leq \frac{1}{2}|\mathcal{H}|$. So we
have
\begin{equation} \sum_{i=1}^k |\mathcal{A}_i| = \sum_{i=1}^k
|\mathcal{A}_i'| + \sum_{i=1}^k |\mathcal{A}_i^*| \leq
|\mathcal{A}'| + k|\mathcal{A}^*| \leq k\left(|\mathcal{A}^*| +
\frac{1}{2}|\mathcal{A}'| \right) \leq k \frac{|\mathcal{H}|}{2}
\nonumber
\end{equation}
and hence, by Lemma~\ref{AMGM},
\begin{align} \prod_{i=1}^k |\mathcal{A}_i| \leq \left(
\frac{1}{k} \sum_{i=1}^k|\mathcal{A}_i| \right)^k \leq \left(
\frac{|\mathcal{H}|}{2} \right)^k.\nonumber
\end{align}

The second part of the theorem is an immediate consequence of
Proposition~\ref{Bergeprop}.~\hfill{$\Box$}

\begin{cor} Conjectures~\ref{weaksum}, \ref{strongsum} and
\ref{prodconj} are true if the bases of $\mathcal{H}$ have a
common element.
\end{cor}
\textbf{Proof.} If the bases of $\mathcal{H}$ have a common
element $x$, then by Corollary~\ref{Bergecor} and
Proposition~\ref{Bergeprop}, $\mathcal{H}\langle x \rangle$ is a
largest star of $\mathcal{H}$ of size $|\mathcal{H}|/2$. By
Theorem~\ref{result2}, the result follows.~\hfill{$\Box$}

\begin{cor} Conjecture~\ref{strongsum} is true if $k = 2$.
\end{cor}
\textbf{Proof.} By Corollary~\ref{Bergecor}, we have
$|\mathcal{S}| \leq |\mathcal{H}|/2$ and hence $2 \leq
\frac{|\mathcal{H}|}{|\mathcal{S}|}$. Now by
Theorem~\ref{result2}, $|\mathcal{A}_1| + |\mathcal{A}_2| \leq
|\mathcal{H}|$. Hence the result.~\hfill{$\Box$}
\\

We now come to our main result, which verifies
Conjectures~\ref{weaksum}, \ref{strongsum} and \ref{prodconj} for
the case when $k \geq n+1$ and $\mathcal{H}$ is compressed with
respect to an element of $[n]$. As remarked in Section~\ref{Her},
an important example of such a hereditary family is one whose
bases have a common element. Other important examples include
$\bigcup_{r = 0}^m{[n] \choose r}$ for any $m \in \{0\} \cup [n]$
(for $m = n$ we get $2^{[n]}$).


\begin{theorem} \label{mainthm} Let $\mathcal{H}$ be a hereditary
sub-family of $2^{[n]}$ that is compressed with respect to an
element $x$ of $[n]$, and let $\mathcal{S} = \mathcal{H}\langle x
\rangle$. Let $k \geq n+1$, and let $\mathcal{A}_1, ...,
\mathcal{A}_k$ be cross-intersecting sub-families of
$\mathcal{H}$. Then
\[\sum_{i = 1}^k |\mathcal{A}_i| \leq k|\mathcal{S}| \quad
\mbox{ and } \quad \prod_{i = 1}^k |\mathcal{A}_i| \leq
|\mathcal{S}|^k,\]
and both bounds are attained if $\mathcal{A}_1 = ... =
\mathcal{A}_k = \mathcal{S}$. Moreover:\\
(a) $\sum_{i = 1}^k |\mathcal{A}_i| = k|\mathcal{S}|$ if and only
if either $\mathcal{A}_1 = ... = \mathcal{A}_k = \mathcal{L}$ for
some largest intersecting sub-family $\mathcal{L}$ of
$\mathcal{H}$ or $k = n+1$ and for some $i \in [k]$,
$\mathcal{A}_i = \mathcal{H} = \{\emptyset\} \cup {[n] \choose
1}$ and $\mathcal{A}_j = \emptyset$ for each $j \in [k]
\backslash \{i\}$.\\
(b) $\prod_{i = 1}^k |\mathcal{A}_i| = |\mathcal{S}|^k$ if and
only if $\mathcal{A}_1 = ... = \mathcal{A}_k = \mathcal{L}$ for
some largest intersecting sub-family $\mathcal{L}$ of
$\mathcal{H}$.
\end{theorem}
This generalises Theorem~\ref{Snevily} in the same way that
Conjectures~\ref{weaksum}, \ref{strongsum} and \ref{prodconj}
generalise Conjecture~\ref{Chvatal} (as explained above). We
prove this result in Section~\ref{Proof1}; however, we set up the
necessary tools in the next section.






\section{New properties of the compression operation}
\label{Compsection}
The proof of Theorem~\ref{mainthm} will be based on the
compression technique, which featured in the original proof of
the classical Erd\H os-Ko-Rado Theorem \cite{EKR}.

For a non-empty set $X$ and $x, y \in X$, let $\delta_{x,y}
\colon 2^X \rightarrow 2^X$ be defined by
\[ \delta_{x,y}(A) = \left\{ \begin{array}{ll}
(A \backslash \{y\}) \cup \{x\} & \mbox{if $y \in A$ and $x \notin
A$};\\
A & \mbox{otherwise,}
\end{array} \right. \]
and let $\Delta_{x,y} \colon 2^{2^X} \rightarrow 2^{2^X}$ be the
\emph{compression operation} (see \cite{EKR}) defined by
\[\Delta_{x,y}(\mathcal{A}) = \{\delta_{x,y}(A) \colon A \in
\mathcal{A}, \delta_{x,y}(A) \notin \mathcal{A}\} \cup \{A \in
\mathcal{A} \colon \delta_{x,y}(A) \in \mathcal{A}\}.\]
Note that $|\Delta_{x,y}(\mathcal{A})| = |\mathcal{A}|$. It is
well-known, and easy to check, that $\Delta_{x,y}(\mathcal{A})$
is intersecting if $\mathcal{A}$ is intersecting; \cite{F}
provides a survey on the properties and uses of compression (also
called \emph{shifting}) operations in extremal set theory. We now
establish new properties of compressions for the purpose of
proving Theorem~\ref{mainthm}.

For any family $\mathcal{A}$, let $\mathcal{A}^*$ denote the
sub-family of $\mathcal{A}$ consisting of those sets in
$\mathcal{A}$ that intersect each set in $\mathcal{A}$
(i.e.~$\mathcal{A}^* = \{A \in \mathcal{A} \colon A \cap B \neq
\emptyset \mbox{ for any } B \in \mathcal{A}\}$), and let
$\mathcal{A}' = \mathcal{A} \backslash \mathcal{A}^*$. So
$\mathcal{A}'$ consists of those sets in $\mathcal{A}$ that do
not intersect all the sets in $\mathcal{A}$, and $\mathcal{A}^*$
is an intersecting family.

\begin{lemma}\label{complemma} Let $\mathcal{A}$ be a sub-family
of $2^{[n]}$, and let $\mathcal{B} = \Delta_{i,j}(\mathcal{A})$
for some $i, j \in [n]$, $i \neq j$. Then: \\
(i) if $A \in \mathcal{A}^*$ then $\delta_{i,j}(A) \in
\mathcal{B}^*$; \\
(ii) if $A \in \mathcal{A}^* \backslash \mathcal{B}^*$ then
$\delta_{i,j}(A) \notin \mathcal{A}^*$; \\
(iii) if $B \in \mathcal{B}^*$ then $\delta_{i,j}(B) \in
\mathcal{B}^*$;\\
(iv) $|\mathcal{A}^*| \leq |\mathcal{B}^*|$.
\end{lemma}
\textbf{Proof.} The lemma is obvious if $\mathcal{A}^* =
\emptyset$, so suppose $\mathcal{A}^* \neq \emptyset$. Fix $A \in
\mathcal{A}^*$.

Obviously $\delta_{i,j}(A) \in \mathcal{B}$. Suppose
$\delta_{i,j}(A) \notin \mathcal{B}^*$. Then $\delta_{i,j}(A) \cap
C = \emptyset$ for some $C \in \mathcal{B}$. By definition of
$\mathcal{B}$, $\delta_{i,j}(C)$ is also in $\mathcal{B}$, and
hence both $C$ and $\delta_{i,j}(C)$ are in $\mathcal{A}$. So $A$
intersects both $C$ and $\delta_{i,j}(C)$. From $\delta_{i,j}(A)
\cap C = \emptyset$ and $A \cap C \neq \emptyset$ we get $i \notin
C$, $\delta_{i,j}(A) \neq A$ (so $i \notin A$), $A \cap C =
\{j\}$. But this yields the contradiction that $A \cap
\delta_{i,j}(C) = \emptyset$. Hence (i).

Suppose $A \notin \mathcal{B}^*$. Assume that $\delta_{i,j}(A) \in
\mathcal{A}^*$. Then $A \in \mathcal{B}$ (as both $A$ and
$\delta_{i,j}(A)$ are in $\mathcal{A}$) and $A \cap D = \emptyset$
for some $D \in \mathcal{B}$ (as $A \notin \mathcal{B}^*$). Since
$A$ intersects each set in $\mathcal{A}$, we must have $D =
\delta_{i,j}(E) \neq E$ for some $E \in \mathcal{A}$, $A \cap E =
\{j\}$ and $i \notin A \cup E$. But then $\delta_{i,j}(A) \cap E =
\emptyset$, contradicting $\delta_{i,j}(A) \in \mathcal{A}^*$. So
$\delta_{i,j}(A) \notin \mathcal{A}^*$. Hence (ii).

Suppose $B \in \mathcal{B}^*$. If $\delta_{i,j}(B) = B$ then
obviously $\delta_{i,j}(B) \in \mathcal{B}^*$. Suppose
$\delta_{i,j}(B) \neq B$. Then $B, \delta_{i,j}(B) \in
\mathcal{A}$. Thus, since $B$ intersects every set in
$\mathcal{B}$ and $i \notin B$, $B$ intersects every set in
$\mathcal{A}$, and hence $B \in \mathcal{A}^*$. By (i),
$\delta_{i,j}(B) \in \mathcal{B}^*$. Hence (iii).

By (i), we can define a function $f : \mathcal{A}^* \rightarrow
\mathcal{B}^*$ by
\[ f(A) =
\begin{cases}\displaystyle
A & \text{if $A \in \mathcal{A}^* \cap \mathcal{B}^*$};\\[2mm]
\delta_{i,j}(A) & \text{if $A \in \mathcal{A}^* \backslash
\mathcal{B}^*$}.
\end{cases} \]
Suppose $A_1, A_2 \in \mathcal{A}^*$ such that $f(A_1) = f(A_2)$.
Suppose $A_1 \in \mathcal{A}^* \cap \mathcal{B}^*$ and $A_2 \in
\mathcal{A}^* \backslash \mathcal{B}^*$; then we have
$\delta_{i,j}(A_2) = f(A_2) = f(A_1) = A_1 \in \mathcal{A}^*$,
which is a contradiction because $\delta_{i,j}(A_2) \notin
\mathcal{A}^*$ by (ii). Similarly, we cannot have $A_2 \in
\mathcal{A}^* \cap \mathcal{B}^*$ and $A_1 \in \mathcal{A}^*
\backslash \mathcal{B}^*$. If $A_1, A_2 \in \mathcal{A}^* \cap
\mathcal{B}^*$ then we have $A_1 = f(A_1) = f(A_2) = A_2$.
Finally, suppose $A_1, A_2 \in \mathcal{A}^* \backslash
\mathcal{B}^*$. Then we have $\delta_{i,j}(A_1) = f(A_1) = f(A_2)
= \delta_{i,j}(A_2)$ and, by (ii), $\delta_{i,j}(A_1) \neq A_1$
and $\delta_{i,j}(A_2) \neq A_2$. So $A_1 =
\delta_{j,i}(\delta_{i,j}(A_1)) = \delta_{j,i}(\delta_{i,j}(A_2))
= A_2$. Therefore, no two distinct sets in $\mathcal{A}^*$ are
mapped by $f$ to the same set in $\mathcal{B}^*$ (i.e.~$f$ is
injective). Hence (iv).
%

\section{Proof of Theorem~\ref{mainthm}} \label{Proof1}

We now prove Theorem~\ref{mainthm}. We adopt the strategy
introduced in \cite{Borg4, Borg5} and also adopted in
\cite{Borg6, BL2, Borg2}, which mainly is to determine, for the
family $\mathcal{F}$ under consideration, the largest real number
$c \leq l/|\mathcal{F}|$ such that $|\mathcal{A}^*| + c
|\mathcal{A}'| \leq l$ for any sub-family $\mathcal{A}$ of
$\mathcal{F}$, where $l$ is the size of a largest intersecting
sub-family of $\mathcal{F}$, and $\mathcal{A}^*$ and
$\mathcal{A}'$ are as defined in Section~\ref{Compsection}; see
\cite{Borg2} for a detailed general explanation.


\begin{theorem} \label{mainthm2} Let $\mathcal{H}$ be a hereditary
sub-family of $2^{[n]}$ that is compressed with respect to an
element $x$ of $[n]$, and let $\mathcal{A}$ be a sub-family of
$\mathcal{H}$. Then
\[|\mathcal{A}^*| + \frac{1}{n+1}|\mathcal{A}'| \leq
|\mathcal{H} \langle x \rangle|,\]
and if $\mathcal{A}' \neq \emptyset$, then equality holds if and
only if $\mathcal{A} = \mathcal{H} = \{\emptyset\} \cup {[n]
\choose 1}$.
\end{theorem}
\textbf{Proof.} Since $\mathcal{H}$ is compressed with respect to
$x$, we have $x \in U(\mathcal{H})$ and hence $\mathcal{H} \langle
x \rangle \neq \emptyset$. The result is trivial if $n=1$, so we
consider 
$n \geq 2$ and proceed by induction on $n$.

We may assume that $x = 1$. Let $\mathcal{L} = \mathcal{H}
\langle 1 \rangle$. Let $\mathcal{B} =
\Delta_{1,n}(\mathcal{A})$. Given that $\mathcal{H}$ is
compressed with respect to $1$, we have $\mathcal{B} \subseteq
\mathcal{H}$.
%
%
%
Define
\begin{align} \mathcal{B}_1 &= \{B \in \mathcal{B} \colon n \in
B\},\nonumber \\
\mathcal{B}_2 &= \{B \backslash \{n\} \colon B \in
\mathcal{B}_1\}, \nonumber \\
\mathcal{B}_3 &= \mathcal{B} \backslash \mathcal{B}_1 = \{B \in
\mathcal{B} \colon n \notin B\}. \nonumber
\end{align}
Define $\mathcal{L}_1, \mathcal{L}_2, \mathcal{L}_3$ and
$\mathcal{H}_1, \mathcal{H}_2, \mathcal{H}_3$ similarly. So
$\mathcal{B}_2, \mathcal{L}_2 \subseteq \mathcal{H}_2 \subseteq
2^{[n-1]}$ and $\mathcal{B}_3, \mathcal{L}_3 \subseteq
\mathcal{H}_3 \subseteq 2^{[n-1]}$. Also note that the properties
of $\mathcal{H}$ are inherited by $\mathcal{H}_3$, that is,
$\mathcal{H}_3$ is hereditary and compressed with respect to $1$;
the same holds for $\mathcal{H}_2$ unless $U(\mathcal{H}_2) =
\emptyset$ (in which case $\mathcal{H}_1$ is either $\emptyset$
or $\{\{n\}\}$). Define
\begin{align} \mathcal{C}_1 &= \{C \in \mathcal{B}_1 \colon 1 \in
C, \mbox{ $B \cap C = \{n\}$ for some } B \in
\mathcal{B}_1\}, \nonumber \\
\mathcal{C}_2 &= \{C \backslash \{n\} \colon C \in \mathcal{C}_1,
\, C \backslash \{n\} \notin \mathcal{B}_3\},
\nonumber \\
\mathcal{D} &= \mathcal{B}_2 \backslash \mathcal{C}_2, \nonumber
\\
\mathcal{E} &= \mathcal{B}_3 \cup \mathcal{C}_2. \nonumber
\end{align}
Obviously $\mathcal{C}_2 \subseteq \mathcal{B}_2$ and $\mathcal{D}
\subseteq \mathcal{H}_2$. Given that $\mathcal{H}$ is hereditary,
we clearly have $\mathcal{C}_2 \subseteq \mathcal{H}_3$; so
$\mathcal{E} \subseteq \mathcal{H}_3$. Note that $\mathcal{L}_2 =
\mathcal{H}_2 \langle 1 \rangle$ and $\mathcal{L}_3 =
\mathcal{H}_3 \langle 1 \rangle$.
Therefore, by the inductive hypothesis, we have $|\mathcal{E}^*| +
\frac{1}{n}|\mathcal{E}'| \leq |\mathcal{L}_3|$, and if
$U(\mathcal{H}_2) \neq \emptyset$, then $|\mathcal{D}^*| +
\frac{1}{n}|\mathcal{D}'| \leq |\mathcal{L}_2|$.

By definition of $\mathcal{C}_2$, we have $\mathcal{B}_3 \cap
\mathcal{C}_2 = \emptyset$ and hence $|\mathcal{E}| =
|\mathcal{B}_3| + |\mathcal{C}_2|$. Since $\mathcal{C}_2 \subseteq
\mathcal{B}_2$, we have $|\mathcal{D}| = |\mathcal{B}_2| -
|\mathcal{C}_2|$. So $|\mathcal{D}| + |\mathcal{E}| =
|\mathcal{B}_2| + |\mathcal{B}_3|$ and hence, since $|\mathcal{D}|
+ |\mathcal{E}| = |\mathcal{D}^*| + |\mathcal{D}'| +
|\mathcal{E}^*| + |\mathcal{E}'|$ and $|\mathcal{B}_2| +
|\mathcal{B}_3| = |\mathcal{B}| = |\mathcal{B}^*| +
|\mathcal{B}'|$,
\begin{equation} |\mathcal{D}^*| + |\mathcal{E}^*| +
|\mathcal{D}'| + |\mathcal{E}'| = |\mathcal{B}^*| +
|\mathcal{B}'|. \label{eq}
\end{equation}

We now come to our main step, which is to show that
$|\mathcal{B}^*| \leq |\mathcal{D}^*| + |\mathcal{E}^*|$. So
suppose $\mathcal{B}^*$ contains a set $B$.

First, suppose $n \notin B$. Then clearly $B$ intersects all sets
in $\mathcal{B}_2 \cup \mathcal{B}_3$ and hence, since
$\mathcal{C}_2 \subseteq \mathcal{B}_2$, we have $B \in
\mathcal{E}^*$. Also, $B \notin \mathcal{C}_2$ since $B \in
\mathcal{B}_3$. In brief, we have
\begin{equation} n \notin B \in \mathcal{B}^* \quad
\Rightarrow \quad B \in (\mathcal{E}^* \backslash \mathcal{C}_2)
\cap \mathcal{B}^*. \label{part1}
\end{equation}
Now suppose $n \in B$, that is, $B \in \mathcal{B}_1$. Let $B^- =
B \backslash \{n\}$. Clearly $B^-$ intersects all sets in
$\mathcal{B}_3$. If $B^- \in \mathcal{C}_2$ then, since all sets
in $\mathcal{C}_2$ contain $1$, $B^-$ also intersects each set in
$\mathcal{C}_2$, meaning that $B^- \in \mathcal{E}^*$. In brief,
we have
\begin{equation} n \in B \in \mathcal{B}^*, \, B \backslash
\{n\} \in \mathcal{C}_2 \quad \Rightarrow \quad B \backslash \{n\}
\in \mathcal{E}^* \cap \mathcal{C}_2. \label{part2}
\end{equation}
Suppose $B^- \notin \mathcal{C}_2$. Then $B^- \in \mathcal{D}$. So
suppose $B^- \notin \mathcal{D}^*$. Then $B^- \cap D = \emptyset$
for some $D \in \mathcal{D}$ and hence, setting $D^+ = D \cup
\{n\}$, we have $B \cap D^+ = \{n\}$ and $D^+ \in \mathcal{B}_1$.
Since $B \cap D = \emptyset$, $D$ cannot be in $\mathcal{B}_3$.
Thus we must have $1 \notin D^+$ because otherwise we get $D^+
\in \mathcal{C}_1$ and hence $D \in \mathcal{C}_2$ (contradicting
$D \in \mathcal{D}$). It follows that we must also have $1 \in B$
because otherwise we get $\delta_{1,n}(B) \cap D^+ = \emptyset$,
contradicting Lemma~\ref{complemma}(iii). So $B \in
\mathcal{C}_1$. Thus, since $B^- \in \mathcal{D}$ implies $B^-
\notin \mathcal{C}_2$, $B^-$ must be in $\mathcal{B}_3$ and hence
in $\mathcal{B}$. Since $B \cap D^+ = \{n\}$, we have $B^- \cap
D^+ = \emptyset$ and hence $B^- \notin \mathcal{B}^*$. However,
since $B^-$ intersects all sets in $\mathcal{B}_3$ and $1 \in B^-
\cap C$ for any $C \in \mathcal{C}_2$, we have $B^- \in
\mathcal{E}^*$. So we have just shown that
\begin{equation} n \in B \in \mathcal{B}^*, \, B \backslash
\{n\} \notin \mathcal{C}_2, \, B \backslash \{n\} \notin
\mathcal{D}^* \quad \Rightarrow \quad B \backslash \{n\} \in
\mathcal{E}^* \backslash (\mathcal{C}_2 \cup \mathcal{B}^*).
\label{part3}
\end{equation}

Define
\begin{align} \mathcal{F}_1 &= \{F \in \mathcal{B}^* \colon n
\notin F\}, \nonumber \\
\mathcal{F}_2 &= \{F \in \mathcal{B}^* \colon n \in F, F
\backslash \{n\} \in \mathcal{C}_2\}, \nonumber \\
\mathcal{F}_3 &= \{F \in \mathcal{B}^* \colon n \in F, F
\backslash \{n\} \notin \mathcal{C}_2, F \backslash \{n\} \notin
\mathcal{D}^*\}, \nonumber \\
\mathcal{F}_4 &= \{F \in \mathcal{B}^* \colon n \in F, F
\backslash \{n\} \notin \mathcal{C}_2, F \backslash \{n\} \in
\mathcal{D}^*\}. \nonumber
\end{align}
Clearly $|\mathcal{B}^*| = |\mathcal{F}_1| + |\mathcal{F}_2| +
|\mathcal{F}_3| + |\mathcal{F}_4|$ and $|\mathcal{F}_4| \leq
|\mathcal{D}^*|$. Also, by (\ref{part1}) - (\ref{part3}), we have
$|\mathcal{F}_1| \leq |(\mathcal{E}^* \backslash \mathcal{C}_2)
\cap \mathcal{B}^*|$, $|\mathcal{F}_2| \leq |\mathcal{E}^* \cap
\mathcal{C}_2|$ and $|\mathcal{F}_3| \leq |\mathcal{E}^*
\backslash (\mathcal{C}_2 \cup \mathcal{B}^*)|$. Thus, since
$(\mathcal{E}^* \backslash \mathcal{C}_2) \cap \mathcal{B}^*$,
$\mathcal{E}^* \cap \mathcal{C}_2$ and $\mathcal{E}^* \backslash
(\mathcal{C}_2 \cup \mathcal{B}^*)$ are disjoint sub-families of
$\mathcal{E}^*$, we obtain $|\mathcal{F}_1| + |\mathcal{F}_2| +
|\mathcal{F}_3| \leq |\mathcal{E}^*|$. So $|\mathcal{B}^*| \leq
|\mathcal{D}^*| + |\mathcal{E}^*|$ as required.

We now know that $|\mathcal{D}^*| + |\mathcal{E}^*| =
|\mathcal{B}^*| + p$ for some integer $p \geq 0$. By (\ref{eq}),
we therefore have $|\mathcal{D}'| + |\mathcal{E}'| =
(|\mathcal{B}^*| + |\mathcal{B}'|) - (|\mathcal{B}^*| + p) =
|\mathcal{B}'| - p$.

At this point, we need to divide the problem into two cases.

\textit{Case 1}: $U(\mathcal{H}_2) \neq \emptyset$. So
$|\mathcal{D}^*| + \frac{1}{n}|\mathcal{D}'| \leq
|\mathcal{L}_2|$. Since we earlier obtained $|\mathcal{E}^*| +
\frac{1}{n}|\mathcal{E}'| \leq |\mathcal{L}_3|$,
\begin{equation} |\mathcal{D}^*| + |\mathcal{E}^*| +
\frac{1}{n}(|\mathcal{D}'| + |\mathcal{E}'|) \leq |\mathcal{L}_2|
+ |\mathcal{L}_3| = |\mathcal{L}_1| + |\mathcal{L}_3| =
|\mathcal{L}|.  \nonumber
\end{equation}
We now have
\begin{equation} |\mathcal{B}^*| + \frac{1}{n}|\mathcal{B}'|
\leq |\mathcal{B}^*| + p + \frac{1}{n}(|\mathcal{B}'| - p) =
|\mathcal{D}^*| + |\mathcal{E}^*| + \frac{1}{n}(|\mathcal{D}'| +
|\mathcal{E}'|) \leq |\mathcal{L}|. \nonumber
\end{equation}
Since $|\mathcal{A}^*| + |\mathcal{A}'| = |\mathcal{A}| =
|\mathcal{B}| = |\mathcal{B}^*| + |\mathcal{B}'|$,
Lemma~\ref{complemma}(iv) gives us $|\mathcal{A}^*| +
\frac{1}{n}|\mathcal{A}'| \leq |\mathcal{B}^*| +
\frac{1}{n}|\mathcal{B}'|$. So $|\mathcal{A}^*| +
\frac{1}{n+1}|\mathcal{A}'| \leq |\mathcal{L}|$, and the
inequality is strict if $\mathcal{A}' \neq \emptyset$.

\textit{Case 2}: $U(\mathcal{H}_2) = \emptyset$. So
$\mathcal{H}_1$ is either $\emptyset$ or $\{\{n\}\}$. If
$\mathcal{H}_1 = \emptyset$, then $\mathcal{H} \subseteq
2^{[n-1]}$ and hence the result follows by inductive hypothesis.
Now suppose $\mathcal{H}_1 = \{\{n\}\}$. Then, since
$\mathcal{B}_1 \subseteq \mathcal{H}_1$, we have $\mathcal{C}_1 =
\mathcal{C}_2 = \emptyset$, which gives $\mathcal{D} =
\mathcal{B}_2 \subseteq \{\emptyset\}$ and $\mathcal{E} =
\mathcal{B}_3$.
If $\mathcal{D} = \emptyset$ then the argument in Case 1 gives us
the result.

Suppose $\mathcal{D} = \{\emptyset\}$. Since $\mathcal{D} =
\mathcal{B}_2$, we have $\mathcal{B}_1 = \{\{n\}\}$ and hence
$\{n\} \in \mathcal{B}$. By definition of $\mathcal{B}$, $\{1\}$
is also in $\mathcal{B}$. Therefore $\mathcal{B} \neq
\mathcal{B}^*$; moreover, since there is no set in $\mathcal{B}
\backslash \{\{n\}\}$ intersecting $\{n\}$, $\mathcal{B} =
\mathcal{B}'$. Now consider $\mathcal{H}_3$. From $\{1\} \in
\mathcal{B} \subseteq \mathcal{H}$ we get $\{1\} \in
\mathcal{H}_3$ and hence $\mathcal{H}_3 \neq \{\emptyset\}$.
Since $\mathcal{H}$ is hereditary, we have $\{\emptyset\} \in
\mathcal{H}_3$, meaning that ${\mathcal{H}_3}^* = \emptyset$ and
$\mathcal{H}_3 = {\mathcal{H}_3}'$. It follows by the inductive
hypothesis that $\frac{1}{n}|\mathcal{H}_3| \leq |\mathcal{L}_3|$
(and hence $n|\mathcal{L}_3| - |\mathcal{H}_3| \geq 0$) and that
equality holds only if $\mathcal{H}_3 = \{\emptyset\} \cup {[n-1]
\choose 1}$. So we have

\begin{align} |\mathcal{L}_3| - \left(|\mathcal{B}^*| +
\frac{1}{n+1}|\mathcal{B}'|\right) &= |\mathcal{L}_3| -
\frac{1}{n+1}|\mathcal{B}| = |\mathcal{L}_3| -
\frac{1}{n+1}(|\mathcal{B}_1| + |\mathcal{B}_3|) \nonumber \\
&\geq |\mathcal{L}_3| - \frac{1}{n+1}(1 + |\mathcal{H}_3|) =
\frac{1}{n+1}((n+1)|\mathcal{L}_3| - 1 - |\mathcal{H}_3|)
\nonumber \\
&= \frac{1}{n+1}(n|\mathcal{L}_3| - |\mathcal{H}_3| +
|\mathcal{L}_3| - 1) \geq \frac{1}{n+1}(|\mathcal{L}_3| - 1) \geq
0, \label{final calc}
\end{align}
where the last inequality follows from the fact that $\{1\} \in
\mathcal{B} \subseteq \mathcal{H}$ and hence $\{1\} \in
\mathcal{L}_3$. So $|\mathcal{B}^*| + \frac{1}{n+1}|\mathcal{B}'|
\leq |\mathcal{L}|$. As in Case 1, we have $|\mathcal{A}^*| +
\frac{1}{n+1}|\mathcal{A}'| \leq |\mathcal{L}|$ by
Lemma~\ref{complemma}.  Suppose equality holds.  Then
$|\mathcal{B}^*| + \frac{1}{n+1}|\mathcal{B}'| = |\mathcal{L}|$.
By the calculation in (\ref{final calc}), we must therefore have
$\frac{1}{n}|\mathcal{H}_3| = \mathcal{L}_3$, implying that
$\mathcal{H}_3 = \{\emptyset\} \cup {[n-1] \choose 1}$, and also
$|\mathcal{B}_3| = |\mathcal{H}_3|$, implying that $\mathcal{B}_3
= \mathcal{H}_3$. Since $\mathcal{B}_1 = \mathcal{H}_1 =
\{\{n\}\}$, $\mathcal{B} = \mathcal{H} = \{\emptyset\} \cup {[n]
\choose 1}$. It clearly follows that $\mathcal{A} = \mathcal{B}$.

Finally, if $\mathcal{A} = \mathcal{H} = \{\emptyset\} \cup {[n]
\choose 1}$, then $\mathcal{A}^* = \emptyset$, $\mathcal{A}' =
\mathcal{A}$, $\mathcal{L} = \{\{1\}\}$, and hence
$|\mathcal{A}^*| + \frac{1}{n+1}|\mathcal{A}'| = 1 =
|\mathcal{L}|$.~\hfill{$\Box$}\\


Now for any non-empty family $\mathcal{F}$, let $l(\mathcal{F})$
be the size of a largest intersecting sub-family of
$\mathcal{F}$, and let $\beta(\mathcal{F})$ be the largest real
number $c \leq l(\mathcal{F})/|\mathcal{F}|$ such that
$|\mathcal{A}^*| + c |\mathcal{A}'| \leq l(\mathcal{F})$ for any
sub-family $\mathcal{A}$ of $\mathcal{F}$.
\\
%
%
\\
\textbf{Proof of Theorem~\ref{mainthm}.} For any intersecting
family $\mathcal{A} \neq \{\emptyset\}$, $\mathcal{A}^* =
\mathcal{A}$ and $\mathcal{A}' = \emptyset$. Thus, by
Theorem~\ref{mainthm2}, $\mathcal{H}\langle x \rangle$ is a
largest intersecting sub-family of $\mathcal{H}$ and hence
$l(\mathcal{H}) = |\mathcal{S}|$. Also by Theorem~\ref{mainthm2},
$\beta(\mathcal{H}) \geq \frac{1}{n+1}$. So we have $k \geq n+1
\geq \frac{1}{\beta(\mathcal{H})}$.

As in the proof of Theorem~\ref{result2}, let $\mathcal{A} =
\bigcup_{i = 1}^k \mathcal{A}_i$; so $\mathcal{A}^* =
\bigcup_{i=1}^k \mathcal{A}_{i}^*$, $\mathcal{A}' = \bigcup_{i =
1}^k \mathcal{A}_i'$ and $|\mathcal{A}'| = \sum_{i=1}^k
|\mathcal{A}_i'|$.
So we have
\begin{equation} \sum_{i=1}^k |\mathcal{A}_i| = \sum_{i=1}^k
|\mathcal{A}_i'| + \sum_{i=1}^k |\mathcal{A}_i^*| \leq
|\mathcal{A}'| + k|\mathcal{A}^*| \leq k\left(|\mathcal{A}^*| +
\beta(\mathcal{H})|\mathcal{A}'| \right) \leq kl(\mathcal{H}) =
k|\mathcal{S}| \label{sumcalc}
\end{equation}
and hence, by Lemma~\ref{AMGM},
\begin{align} \prod_{i=1}^k |\mathcal{A}_i| \leq \left(
\frac{1}{k} \sum_{i=1}^k|\mathcal{A}_i| \right)^k \leq
|\mathcal{S}|^k \label{prodcalc}.
\end{align}

We now prove (a) and (b). It is trivial that the conditions in (a)
and (b) are sufficient, so it remains to prove that they are also
necessary.

Consider first $k > n+1$. Then $k > \frac{1}{\beta(\mathcal{H})}$.
From (\ref{sumcalc}) we see that $\sum_{i=1}^k |\mathcal{A}_i| =
k|\mathcal{S}|$ only if $|\mathcal{A}'| = 0$ and
$|\mathcal{A}_1^*| = ... = |\mathcal{A}_k^*| = |\mathcal{A}^*| =
|\mathcal{S}|$. So $\sum_{i=1}^k |\mathcal{A}_i| = k|\mathcal{S}|$
only if $\mathcal{A}$ is a largest intersecting sub-family of
$\mathcal{H}$ and $\mathcal{A}_1 = ... = \mathcal{A}_k =
\mathcal{A}$. It follows from (\ref{prodcalc}) that $\prod_{i=1}^k
|\mathcal{A}_i| = |\mathcal{S}|^k$ only if $\mathcal{A}$ is a
largest intersecting sub-family of $\mathcal{H}$ and
$\mathcal{A}_1 = ... = \mathcal{A}_k = \mathcal{A}$.

Now consider $k = n+1$. If we still have $k >
\frac{1}{\beta(\mathcal{H})}$, then we arrive at the same
conclusion as in the previous case $k > n+1$. So suppose $k =
\frac{1}{\beta(\mathcal{H})}$. Then $\beta(\mathcal{H}) =
\frac{1}{n+1}$.

Suppose $\mathcal{H} \neq \{\emptyset\} \cup {[n] \choose 1}$.
Let $d = \frac{l(\mathcal{H})}{|\mathcal{H}|}$. Since $x \in
U(\mathcal{H})$, $\mathcal{S} \neq \emptyset$. Thus, since
$|\mathcal{H}| = \left|\{\emptyset\} \cup \bigcup_{i=1}^n
\mathcal{H}\langle i \rangle \right|$, we get $|\mathcal{H}| \leq
1+n|\mathcal{S}|$, and equality holds only if
$|\mathcal{H}\langle i \rangle| = |\mathcal{S}|$ for all $i \in
[n]$. So $|\mathcal{H}| \leq (n+1)|\mathcal{S}|$, and equality
holds only if $|\mathcal{S}| = 1$ and $|\mathcal{H}\langle i
\rangle| = |\mathcal{S}|$ for all $i \in [n]$. If $i \in [n]$ and
$A \in \mathcal{H}\langle i \rangle$, then, since $\mathcal{H}$ is
hereditary, all subsets of $A$ containing $i$ are also in
$\mathcal{H}\langle i \rangle$; thus, if $|\mathcal{H}\langle i
\rangle| = 1$, then $\mathcal{H}\langle i \rangle$ must be
$\{i\}$. Therefore, if $|\mathcal{H}| = (n+1)|\mathcal{S}|$, then
$\mathcal{H}\langle i \rangle = \{i\}$ for all $i \in [n]$, but
this gives the contradiction that $\mathcal{H} = \{\emptyset\}
\cup {[n] \choose 1}$. So $|\mathcal{H}| < (n+1)|\mathcal{S}|$
and hence $\frac{1}{n+1} < \frac{|\mathcal{S}|}{|\mathcal{H}|} =
d$. Now let $\mathcal{A} \subseteq \mathcal{H}$. If $\mathcal{A}'
= \emptyset$, then obviously $|\mathcal{A}^*| + d|\mathcal{A}'|
\leq l(\mathcal{H})$. If $\mathcal{A}' \neq \emptyset$, then
$|\mathcal{A}^*| + \frac{1}{n+1}|\mathcal{A}'| < l(\mathcal{H})$
by Theorem~\ref{mainthm2}. Thus, if $c$ is the largest real
number such that $c \leq d$ and $|\mathcal{A}^*| + c
|\mathcal{A}'| \leq l(\mathcal{H})$ for any $\mathcal{A}
\subseteq \mathcal{H}$, then $c > \frac{1}{n+1}$, which is a
contradiction since $\beta(\mathcal{H}) = \frac{1}{n+1}$.


We have therefore shown that $\mathcal{H}$ must consist of the
sets $\emptyset, \{1\}, \{2\}, ..., \{n\}$. It follows by the
cross-intersection condition that we have the following:\\
- If one of the families $\mathcal{A}_1, ..., \mathcal{A}_k$
consists of only one set $A$ and $A \neq \emptyset$, then each of
the others either consists of $A$ only or is empty.\\
- If one of the families $\mathcal{A}_1, ..., \mathcal{A}_k$
either has more than one set or has the set $\emptyset$, then the
others must be empty. \\
These have the following immediate implications:\\
- If $\sum_{i = 1}^k |\mathcal{A}_i| = k|\mathcal{S}|$, then
$\sum_{i = 1}^k |\mathcal{A}_i| = n+1$ (since $\mathcal{S} =
\{\{x\}\}$ and $k = n+1$) and hence either $\mathcal{A}_1 = ... =
\mathcal{A}_k = \{\{y\}\}$ for some $y \in [n]$, or for some $i
\in [n]$, $\mathcal{A}_i = \mathcal{H} = \{\emptyset\} \cup {[n]
\choose 1}$ and $\mathcal{A}_j = \emptyset$ for each $j \in [n]
\backslash \{i\}$.\\
- If $\prod_{i = 1}^k |\mathcal{A}_i| = |\mathcal{S}|^k$, then
$\sum_{i = 1}^k |\mathcal{A}_i| = 1$ and hence $\mathcal{A}_1 =
... = \mathcal{A}_k = \{\{y\}\}$ for some $y \in [n]$.\\
Note that for any $y \in [n]$, $\{\{y\}\}$ is a largest
intersecting sub-family of $\mathcal{H} = \{\emptyset\} \cup {[n]
\choose 1}$.~\hfill{$\Box$}


\end{document}